\documentclass[11pt]{amsart}
\usepackage{epsf}

\newcommand{\epr}{\qed}
\newcommand{\A}{{\mathfrak A}}
\newcommand{\g}{{\mathfrak g}}
\newcommand{\h}{{\mathfrak h}}
\newcommand{\gl}{{\mathfrak{gl}}}

\def\matho#1{\mathop{\mathrm{#1}}}
\def\wt{\widetilde}
\newcommand*{\ato}[2]{{\genfrac{}{}{0pt}{}{#1}{#2}}}
\def\suml{\sum\limits}
\DeclareMathSizes{11.1}{10}{8}{6}

\newcommand\nfrac[2]
{\dfrac{\raisebox{-1pt}{\fontsize{11.1}{10pt}\selectfont$#1$}}
       {\raisebox{ 1pt}{\fontsize{11.1}{10pt}\selectfont$#2$}}}

\newcommand{\Alt}{\matho{Alt}}
\newcommand{\Ker}{\matho{Ker}}
\newcommand{\fin}{\matho{fin}}
\newcommand{\Jac}{\matho{Jac}}
\newcommand{\Res}{\matho{Res}}
\newcommand{\Der}{\matho{Der}}
\newcommand{\Lie}{\matho{Lie}}
\newcommand{\Inn}{\matho{Inn}}
\newcommand{\Ass}{\matho{Ass}}
\newcommand{\codim}{\matho{codim}}
\newcommand{\Cycle}{\matho{Cycle}}
\newcommand{\Altl}{\Alt\limits}
\newcommand{\Tr}{\matho{Tr}}
\newcommand{\pr}{\matho{pr}}
\newcommand{\ad}{\matho{ad}}
\newcommand{\Dif}{\mathrm{Dif}}
\newcommand{\Poiss}{\matho{Poiss}}

\newcommand{\C}{\mathbb C}
\newcommand{\R}{\mathbb R}
\newcommand{\Z}{\mathbb Z}
\newcommand{\D}{\mathcal D}

\newtheorem*{theorem}{Theorem}
\newtheorem*{lemma}{Lemma}
\newtheorem*{corollary}{Corollary}
\newtheorem*{conjecture}{Conjecture}
\newtheorem*{mainconjecture}{Main Conjecture}
\newtheorem*{question}{Question}

\theoremstyle{remark}
\newtheorem*{remark}{Remark}
\newtheorem*{example}{Example}

\author{Boris Shoikhet}
\title[Lifting formulas]%
{Cohomology of the Lie algebras of differential operators: lifting
formulas}
\date{1997}
\address{IUM, 11 Bol'shoj Vlas'evskij per.,
Moscow 121002, Russia}
\email{borya@main.mccme.rssi.ru}

\begin{document}
\maketitle
 \sloppy

\def\nfracp#1#2{\nfrac{\partial#1}{\partial#2}}

\section*{Introduction}

The lifting formulas are formulas for the cocycles of the Lie
algebra $\A$,
constructed from an associative algebra~$\A$ with
the bracket $[a,b]=a\cdot b-b\cdot a$. These cocycles are
constructed with the following data:
\renewcommand{\theenumi}{\arabic{enumi})}
\renewcommand{\labelenumi}{\textup{\theenumi}}
\begin{enumerate}
\item a trace $\Tr\colon\A\to\C$ on the associative algebra~$\A$
\item a set $\D=\{D_1,D_2,\dots\}$ of the (exterior) derivations
of the associative algebra~$\A$
\end{enumerate}
satisfying the following conditions:
\renewcommand{\theenumi}{(\roman{enumi})}
\begin{enumerate}
\item $\Tr(D_iA)=0$ for any $A\in\A$ and any $D_i\in\D$
\item $[D_i,D_j]=\ad(Q_{ij})$ --- inner derivation
($Q_{ij}\in\A$) for any $D_i,D_j\in\D$
\item $\Altl_{i,j,k}D_k(Q_{ij})=0$ for all $i$, $j$, $k$.
\end{enumerate}

The main example of such a situation is the Lie algebra
$\Psi\Dif_n(S^1)$ of the formal pseudodifferential operators on
$(S^1)^n$ (see \cite{A}). The trace $\Tr$ in this example is the
``noncommutative residue'', $\Tr(D)=$ the coefficient of the
term $x_1^{-1}\cdot x_2^{-1}\cdot\ldots\cdot x_n^{-1}\cdot
\partial_1^{-1}\ldots\cdot\nobreak
\partial_n^{-1}$ of $D\in\Psi\Dif_n(S^1)$
(in any coordinate system). It is easy to verify that
$\Tr[D_1,D_2]=0$ for any $D_1,D_2\in\Psi\Dif_n(S^1)$. Furthermore,
$\ln x_i$ ($i=1,\dots,n)$ are (exterior) derivations of
$\Psi\Dif_n(S^1)$ with respect to the adjoint action; the
symmetry between the operators~$x$ and $\nfrac d{dx}$ allows us to
define exterior derivations $\ln\partial_i$ ($i=1,\dots,n$).
(Actually, one can define $\ln D$\ \ ($D\in\Psi\Dif_n(S^1)$) in much
greater generality --- see Sect.\ 2.4 and 2.8).

We prove in \S1 that the noncommutative residue~$\Tr$ on the
associative algebra $\Psi\Dif_n(S^1)$ and the set of $2n$
derivations $\{\ln x_i,\ln\partial_i;i=1,\dots,n\}$ satisfy
conditions (i)--(iii) above, and this is our main (and in some
sense the unique) example.

In the case of the one derivation~$D$ such a construction appeared
in~\cite{KK}, where two $2$-cocycles on the Lie algebra
$\Psi\Dif_1(S^1)=\Psi\Dif(S^1)$ were constructed:
\begin{gather*}
\Psi^{(1)}(A_1,A_2)=\Tr([\ln x,A_1]\cdot A_2)\\
\Psi^{(2)}(A_1,A_2)=\Tr([\ln \partial,A_1]\cdot A_2).
\end{gather*}
Both these cocycles are cohomologous to zero when restricted  to
the  Lie  algebra  $\Dif_1$  of  the  (polynomial)  differential
operators on~$\C^1$; on the other hand, our aim is to  construct
cocycles  on  this Lie algebra. We are able to accomplish this by the
simultaneous application of both $\ln x_i$ and $\ln\partial_i$.

Now the problem is solved only for $n=1$ (\S1) and $n=2$ (\S3).
But the Second Version of the
Main Conjecture (\S4) gives  us an explicit formula for arbitrary~$n$
(when conditions (i)--(iii) above hold). We obtain this
formulas using a some simple trick, and all that remains to prove is
the fact that this formula gives
in fact \emph{cocycles}. We will see in \S3
($n=2$) that this is very difficult problem.

In  this  way we obtain elements in $H^{2n+1}(\Dif_n;\C)$.
We assume that $H^{2n+1}(\Dif_n;\C)=\C$ (but we don't use this fact elsewhere).
There
is a  deformation  of  the  Lie  algebra  $\Poiss_{2n}$  of  the
functions  on~$\R^{2n}$  with  the  Poisson  bracket  to the Lie
algebra $\Dif_n$, and the well-known  result  is  the  vanishing
theorem:
$H^i(\Poiss_{2n};\C)=0$ when
   $0<i\le 2n$. Hence, we study the first nontrivial case.

It is easy to construct a $(2n+1)$-cocycle on
$\Poiss_{2n}$ (or on the Lie algebra $\Poiss_{2n}(S^1)$ of the
functions on the torus $(T^2)^n$ with the Poisson bracket), and
if the lifting problem can be solved, the corresponding
$(2n+1)$-cocycle on $\Psi\Dif_n(S^1)$ is a deformation of this
$(2n+1)$-cocycle on $\Poiss_{2n}(S^1)$. The lifting problem
is equivalent to
the deformation problem~(\S2).

B.\,Feigin explained me that such formulas should exist;
numerous talks with him and with M.\,Kontsevich made me
interested in this subject. I express my deep gratitude to them.
I am also grateful to Jeremy Bem for helping with my English
and to Seva Kordonsky for carefully typing this manuscript.

\section{Cohomology of the Lie algebra $\Dif_1$}

\subsection{}
Let $\Tr=\pr_{x^{-1}\partial^{-1}}$ be a trace on the associative algebra
$\Psi\Dif(S^1)$, $D_1=\ad(\ln\partial)$ and $D_2=\ad(\ln x)$ be
two its exterior derivations.

\begin{lemma}
\begin{enumerate}
\item $[D_1,D_2]=\ad Q$ \emph($Q\in\Psi\Dif(S^1)$\emph) is an inner
derivation of the associative algebra $\Psi\Dif(S^1)$.
\item 
$Q=x^{-1}\partial^{-1}+\nfrac12x^{-2}\partial^{-2}+\nfrac23x^{-3}\partial^{-3}
+
\ldots+\nfrac{(n-1)!}nx^{-n}\partial^{-n}+\dots$
\end{enumerate}
\end{lemma}

\begin{proof}

It is sufficient to show that $[D_1,D_2](x)=[Q,x]$ and
$[D_1,D_2](\partial)=[Q,\partial]$. This is a straightforward
calculation.
\end{proof}

\begin{theorem}
$$
\Psi_3(A_1,A_2,A_3)=\Alt\limits_{A,D}\Tr(D_1A_1\cdot D_2A_2\cdot A_3)+
\Alt\limits_A\Tr(Q\cdot A_1\cdot A_2\cdot A_3)$$
 is a $3$-cocycle
on  the  Lie  algebra  $\Psi\Dif(S^1)$.  This  cocycle  is   not
cohomologous to zero, and it remains noncohomologous to zero after
restriction to the Lie algebra $\Dif_1$.
\end{theorem}

We will prove this theorem in Sect. 1.2--1.4.

\subsection{}

\begin{lemma}

$\Tr(D_iA)=0$ for any $A\in\Psi\Dif(S^1)$ and $i=1,2$.
\end{lemma}

\begin{proof}

It is clear that $[\Psi\Dif(S^1),\Psi\Dif(S^1)]\in\Ker\Tr D_iA$,
but   $\codim([\Psi\Dif(S^1),\Psi\Dif(S^1)])=1$.   Hence  it  is
sufficient to verify the Lemma for $A=x^{-1}\partial^{-1}$.
\end{proof}

\begin{corollary}

\begin{equation}
\Altl_{A,D}\Tr(D_1(D_2A_1\cdot A_2\cdot A_3\cdot A_4))=0
\end{equation}
for any $A_1,A_2,A_3,A_4\in\Psi\Dif(S^1)$.
\end{corollary}

\epr

The alternation here is just  decoration.

\subsection{}

Let us calculate the  l.h.s. of~(1), using the Leibniz rule; we have:
\begin{multline}
\text{l.h.s (1)}=\Altl_{A,D}\Tr(D_2A_1\cdot D_1A_2\cdot A_3\cdot
A_4+D_2A_1\cdot           A_2           \cdot           A_3\cdot
D_1A_4)+\\
+\Altl_A\Tr([Q,A_1]\cdot A_2\cdot A_3\cdot A_4).
\end{multline}
Note   that  $\Altl_{A,D}\Tr(D_2A_1\cdot  A_2\cdot  D_1A_3\cdot
A_4)\equiv0$ in view of  the  main  property  of~$\Tr$  and  the
symmetry conditions,
\begin{equation}
\text{r.h.s. (2)}=-2\Altl_{A,D}\Tr(D_1A_1\cdot D_2A_2\cdot
A_3\cdot A_4)+\Altl_A\Tr([Q,A_1]\cdot A_2\cdot A_3\cdot A_4)
\end{equation}
Let us denote the first summand in (3) by~$\alpha$ and the
second one by~$\beta$

\begin{lemma}
\begin{enumerate}
\item $\alpha=d(-2\Altl_{A,D}\Tr(D_1A_1\cdot D_2A_2\cdot 
A_3))(A_4,A_1,A_2,A_3)$;
\item $\beta=d(-2\Altl_A\Tr(Q\cdot A_1\cdot A_2\cdot A_3))(A_4,A_1,A_2,A_3)$.
\end{enumerate}
\end{lemma}

\begin{proof}

We will use the following expression for the differential in the
cochain complex of the Lie algebra:
\begin{multline}
(d\Psi_n)(A_1,\dots,A_{n+1})=
\nfrac12\cdot\sum_{i=1}^{n+1}(-1)^{i-1}\cdot\\
\cdot(\Psi_n([A_i,A_1],A_2,\dots,\hat
A_i,\dots,A_{n+1})+\Psi_n(
A_1,[A_i,A_2],\dots,\hat
A_i,\dots,A_{n+1})+\dots)
\end{multline}
Let us prove (i):
\begin{multline}
d(\Altl_{A,D}\Tr(D_1A_1\cdot D_2A_2\cdot A_3))(A_4,A_1,A_2,A_3)=\\
=\nfrac12\Altl_\ato{A_1,\dots,A_4}{D_1,D_2}\Tr(D_1[A_4,A_1]\cdot
D_2A_2\cdot A_3+\\
+D_1A_1\cdot D_2[A_4,A_2]\cdot A_3+D_1A_1\cdot
D_2A_2\cdot[A_4,A_3]).
\end{multline}
Then we subtract from r.h.s. of (5)
\begin{multline*}
\gamma=\nfrac12\Altl_{A,D}\Tr([A_4,D_1A_1]\cdot D_2A_2\cdot
A_3+\\
+D_1A_1\cdot[A_4,D_2A_2]\cdot A_3+D_1A_1\cdot
D_2A_2\cdot[A_4,A_3]),
\end{multline*}
(it is clear that $\gamma=0$).

We have:
\begin{gather*}
\text{r.h.s. (5)}=\nfrac12\Altl_\ato{A_1,\dots,A_4}{D_1,D_2}
\Tr([D_1A_4,A_1]\cdot D_2A_2\cdot A_3+D_1A_1\cdot[D_2A_4,A_2]\cdot A_3),\\
\Altl_{A,D}\Tr(D_1A_4\cdot A_1\cdot D_2A_2\cdot A_3-D_1A_1\cdot A_2\cdot 
D_2A_4\cdot A_3)=0
\end{gather*}
in view of the symmetry conditions, and
$$
\text{r.h.s (5)}=\nfrac12\cdot 2\cdot\Altl_{A,D}\Tr(D_1A_1\cdot
D_2A_2\cdot A_3\cdot A_4).
$$

The proof of (ii) is similar.
\end{proof}

\subsection{}

It follows from Corollary 1.2 and Lemma 1.3 that
$$
d(\Altl_{A,D}\Tr(D_1A_1\cdot D_2A_2\cdot A_3)+\Altl_A\Tr(Q\cdot
A_1\cdot A_2\cdot A_3))=0,
$$
and therefore
$$
\Psi_3(A_1,A_2,A_3)=\Altl_{A,D}\Tr(D_1A_1\cdot D_2A_2\cdot
A_3)+\Altl_A\Tr(Q\cdot A_1\cdot A_2\cdot A_3)
$$
is a $3$-cocycle on the Lie algebra $\Psi\Dif(S^1)$.
To complete the proof of Theorem 1.1, it is sufficient to verify
that the restriction of~$\Psi_3$ on $\Dif_1$ is nonhomologous to
zero cocycle. The check is straightforward: $1\wedge
x\wedge\partial$ is a cycle on $\Dif_1$ and $\Psi_3(1\wedge
x\wedge\partial)\ne0$.
\epr

\subsection{}

In \cite{F} B.\, Feigin formulated a conjecture about
$H^*(\Dif_1;\C)$. According to this conjecture,
$$
H^*(\Dif{}_1,\C)=\wedge^*(\Psi_3,\Psi_5,\Psi_7,\dots)\otimes
S^*(c_4,c_6,c_8,\dots)
$$
(the lower index denotes the grading). This conjecture has not been proved
yet. However, our methods allow to construct some
cocycles $\wt\Psi_5$, $\wt\Psi_7$, $\wt\Psi_9$, \dots which are
``candidates'' for $\Psi_5$, $\Psi_7$, $\Psi_9$, \dots;
unfortunately, even the proof of their nontriviality runs into
trouble because there are not explicit formulas for the
cycles. Also, it seems that one can't obtain a proof by passing
to the Hamiltonian limit. Note that no representatives for
$c_4$, $c_6$, $c_8$, \dots have been found.

\begin{theorem}
For $i=2,3,4,\dots$
\begin{multline*}
\wt\Psi_{2i+1}(A_1,\dots,A_{2i+1})=
\Altl_A\Tr(Q\cdot A_1\cdot\ldots\cdot A_{2i+1})+\\
+\Altl_{A,D}\Tr(D_1A_1\cdot D_2A_2\cdot A_3\cdot\ldots
A_{2i+1}+D_1A_1\cdot A_2\cdot A_3\cdot D_2A_4\cdot
A_5\cdot\ldots\cdot A_{2i+1}+\\
+D_1A_1\cdot A_2\cdot A_3\cdot A_4\cdot A_5\cdot D_2A_6\cdot
A_7\cdot\ldots\cdot A_{2i+1}+\ldots+\\
+a_s\cdot D_1A_1\cdot
A_2\cdot\ldots\cdot A_s\cdot D_2A_{s+1}\cdot A_{s+2}\cdot\ldots\cdot A_{2i+1})
\end{multline*}
where $\begin{cases}\text{$s=i+1$,\ \ $a_s=\frac12$, if $i$ is even}\\
\text{$s=i$,\ \ $a_s=1$, if $i$ is odd}\end{cases}$\newline
is a $(2i+1)$-cocycle on the Lie algebra $\Psi\Dif(S^1)$
\emph(or $\Dif_1$\emph).
\end{theorem}

\begin{proof}
It is a direct calculation, similar to the Proof of Theorem 1.1.
\end{proof}

\begin{conjecture}
$\Psi_3$, $\wt\Psi_5$, $\wt\Psi_7$, $\wt\Psi_9$ \dots generate the
exterior algebra in $H^*(\Dif_1,\C)$.
\end{conjecture}

Note, that from this Conjecture follows the statement about
spectral sequence from~\cite{F}. Therefore, if this Conjecture
is true, all we need to prove Feigin's conjecture about
$H^*(\Dif_1,\C)$ is the ``main theorem of the invariant theory''
for the Lie algebra $\gl(\lambda)$ --- see~\cite{F}.

\section{The lifting problem and the deformation problem}

\subsection{}

Let $\A$ be an associative algebra, let us denote by the same
symbol the corresponding Lie algebra with the bracket
$[a,b]=a\cdot b-b\cdot a$. Suppose that there is an invariant
form $\langle,\rangle$ on the Lie algebra~$\A$ (not necessary
nondegenerate).

For any Lie algebra $\g$ there is a canonical map $\theta_{i+1}\colon
H^{i+1}(\g;\C)\to H^i(\g;\g^*)$. Inner derivations $\Inn(\A)$ of
the associative algebra~$\A$ form an ideal in the Lie algebra
$\Der_{\Ass}(\A)$ of all derivations of the associative
algebra~$\A$,\ \ $\Inn(\A)$ is also an ideal in the Lie algebra
$\Der_{\Lie}(\A)$ of the derivations of the Lie algebra~$\A$.
There is a map $\Der_{\Ass}(\A)\to\Der_{\Lie}(\A)$ and a map
$\Der_{\Ass}(\A)/\Inn(\A)\to\Der_{\Lie}(\A)/\Inn(\A)$. The last
space is equal to $H^1(\A,\A)$.
Let $D_1,\dots,D_i\in\Der_{\Ass}(\A)$; we will use the same
symbols for their images in $H^1(\A.\A)$.
There is a map $H^m(\A,\A)\otimes H^n(\A;\A)\to H^{m+n}(\A;\A)$
because $\A$ is an associative algebra.

Let $\alpha_{D_1,\dots,D_i}=\Altl_D(D_1\cdot D_2\cdot\ldots\cdot
D_i)\in H^i(\A;\A)$ or, in the explicit form,
$\alpha_{D_1,\dots,D_i}(A_1,\dots,A_i)=\Altl_{A,D}(D_1(A_1)\cdot
D_2(A_2)\cdot\ldots\cdot D_i(A_i))$. The form $\langle,\rangle$
determines the invariant map $\A\to\A^*$, and we have the element
$\alpha^*_{D_1,\dots,D_i}\in H^i(\A;\A^*)$.

\begin{question}
When does there exist an element $\Psi_{i+1}\in H^{i+1}(\A,\C)$ such
that $\theta_{i+1}(\Psi_{i+1})=\alpha^*_{D_1,\dots,D_i}$. If
such an element exists, it is called \emph{lifting} of
$\alpha^*_{D_1,\dots,D_n}$.
\emph(Note, that $\Psi_{i+1}$, if it exists, is not uniquely determined.\emph)

Furthermore, if $D_i^{(j)}\in H^1(\A;\A)$,\ \ $a_j\in\C$, we have an
element $\suml a_k\alpha^*_{D_1^{(k)},\dots,D_i^{(k)}}$, and
there is the same lifting problem for this element.
\end{question}

\subsection{}

We will be concerned only with the case when there is a trace
$\Tr\colon\A\to\C$ on the associative algebra~$\A$ (i.e.,
$\Tr[A_1,A_2]=0$ for any $A_1,A_2\in\A$) and $\langle
A_1,A_2\rangle=\Tr(A_1\cdot A_2)$. Suppose also, that the
condition~(i) from the Introduction holds: $\Tr(DA)=0$ for all
$A\in\A$ and $D\in\D=\{D_1,D_2,\dots,D_i\}$. Such derivations
form an Lie subalgebra
$\Der_{\Ass}^{\Tr}(\A)\subset\Der_{\Ass}(\A)$; we will denote by
$H^1_{\Tr}(\A;\A)$ the image of this subalgebra in $H^1(\A;\A)$.

\begin{lemma}
Suppose that $\suml_aa_k\cdot D_1^{(k)}\wedge\dots\wedge
D_i^{(k)}$ is a \emph{cycle} in Lie algebra
$\Der_{\Ass}^{\Tr}(\A)$ \emph{(but not only in
$H^1_{\Tr}(\A;\A)$)}. Then
$\Psi_{i+1}(A_1,\dots,A_{i+1})=\Altl_{A,D}\Tr(\suml_ka_kD_1^{(k)}
A_1\cdot\ldots\cdot D_i^{(k)}A_i\cdot A_{i+1})$ is a cocycle,
and the corresponding element in $H^{i+1}(\A;\C)$ is a
\emph{lifting} of
$\suml_ka_k\alpha^*_{D_1^{(k)},\dots,D_i^{(k)}}$.
\end{lemma}

\begin{proof}
It is sufficient to prove that $\Psi_{i+1}$ is a cocycle. As
in the proof of Theorem 1.1 we have to write an expression on
$A_1$, \dots, $A_{i+2}$ which is a coboundary and which a priori
is equal to zero.

For simplicity we suppose, that $D_1\wedge\dots\wedge D_i$ is
a cycle in $\Der_{\Ass}^{\Tr}(\A)$.

Consider the following expression:
\begin{equation}
\Altl_{A,D}\Tr(D_i(\Cycle(D_1A_1\cdot\ldots\cdot
D_{i_1}A_{i_1}\cdot A_i\cdot A_{i+1})\cdot A_{i+2}).
\end{equation}
(Here $\Cycle$ is the sum on all the cyclic permutations; it
contains $i+1$ summands.)

The part of (6) which contains $[D_a,D_b]$ is equal to zero
because $D_1\wedge\dots\wedge D_i$ is a cycle.

The remaining summands are of $3$ types:
\begin{align*}
&\Alt\Tr(D\cdot D\cdot\ldots\cdot D\cdot A\cdot D\cdot D\cdot A),\\
&\Alt\Tr(D\cdot D\cdot\ldots\cdot D\cdot D\cdot A\cdot D\cdot A),\\
\text{and}\ &\Alt\Tr(D\cdot D\cdot\ldots\cdot D\cdot D\cdot D\cdot A\cdot A).\\
\end{align*}
(Here $A$ denotes~$A_a$ and $D$ denotes $D_bA_a$.)
The summands of the first two types are \emph{not}
coboundaries, but the summands of the third type are. Therefore
we have to eliminate the summands of the first two types.

We have:
\begin{multline}
\text{(6)}=\Altl_{A,D}\Tr(\Cycle(D_1A_1\cdot\dots\cdot
D_{i-1}A_{i-1}\cdot D_iA_i\cdot A_{i+1})\cdot A_{i+2})+\\
+\Altl_{A,D}\Tr(\Cycle(D_1A_1\cdot\dots\cdot
D_{i-1}A_{i-1}\cdot A_i\cdot D_iA_{i+1})\cdot A_{i+2})+\\
+\Altl_{A,D}\Tr(\Cycle(D_1A_1\cdot\dots\cdot
D_{i-1}A_{i-1}\cdot A_i\cdot A_{i+1})\cdot D_iA_{i+2})
\end{multline}
The first summand in (7) cancels with the second one because of
the symmetry and the main property of~$\Tr$; the third summand
is
\begin{multline}
i\cdot\Altl_{A,D}\Tr(D_1A_1\cdot\ldots\cdot D_iA_i\cdot A_{i+1}\cdot
A_{i+2})-\\
-\Alt_{A,D}\Tr(D_1A_1\cdot\ldots D_{i-1}A_{i-1}\cdot A_i\cdot D_iA_{i+1}\cdot
A_{i+2}).
\end{multline}
Furthermore let us consider an expression
\begin{equation}
\Altl_{A,D}\Tr(D_i(D_1A_1\cdot\ldots\cdot D_{i-1}A_{i-1}\cdot
A_i\cdot A_{i+1}\cdot A_{i+2})).
\end{equation}
It is clear that
\begin{multline*}
\text{(8)}+\text{(9)}=(i+2)\Altl_{A,D}\Tr(D_1A_1\cdot\ldots\cdot
D_iA_i\cdot A_{i+1}\cdot
A_{i+2})=\\
=(i+2)d(\Altl_{A,D}\Tr(D_1A_1\cdot\ldots\cdot
D_iA_i\cdot A_{i+1}))(A_{i+2},A_1,\dots,A_{i+1}).
\end{multline*}
The terms which contain $[D_a,D_b]$ again cancel because
$D_1\wedge\dots\wedge D_i$ is a cycle.
\end{proof}

\subsection{}

We give here two examples of Lemma 2.2.

\subsubsection{}

\begin{example}
Let $\g=\gl_n\otimes\C[t_1,\dots,t_n,t_1^{-1},\dots,t_n^{-1}]$.
Define $\Tr\colon\g\to\C$ to be the composition
$\Tr=\Res_{t_1,\dots,t_n}\circ\Tr_{\gl_n}$. The associative
algebra~$\g$ has $n$ exterior derivations: $D_1=\nobreak\nfrac d{dt_1}$,
\dots, $D_n=\nfrac d{dt_n}$. Therefore, Lemma 2.2 gives the
$(n+1)$-cocycle on the Lie algebra~$\g$. In particular, in the
case $n=1$ we obtain the standard Kac--Moody $2$-cocycle on the
current algebra.
\end{example}

\subsubsection{}

\begin{example}

(i)
Let $\A$  be an associative algebra, $\h$ an two-sided ideal in~$\A$.
Then~$\A$ acts on the associate algebra~$\h$ via the adjoint
action. Suppose, that the image of~$\A$ in $\Der_{\Ass}(\h)$
lies in $\Der_{\Ass}^{\Tr}(\h)$. Then Lemma 2.2 gives a map
$q\colon H_i(\A;\C)\to H^{i+1}(\h;\C)$.

(ii)
In particular, let $\gl_\infty^{\Jac}$ be an associative algebra
of infinite matrices with a finite number (depending on the
matrix) of nonzero diagonals, and let $\gl_\infty^{\fin}$ be the
ideal in~$\gl_\infty^{\Jac}$ consisting of the matrices with
finite support. Then we have a map
$$
q\colon H_i(\gl_\infty^{\Jac};\C)\to H^{i+1}(\gl_\infty^{\fin};\C).
$$
Recall (\cite{FT,Fu}) that
\begin{gather*}
H^*(\gl_\infty^{\Jac};\C)=S^*(c_2,c_4,c_6,\dots)\\
H^*(\gl_\infty^{\fin};\C)=\Lambda^*(\xi_1,\xi_3,\xi_5,\dots)
\end{gather*}
This construction has an evident generalization: let $A$ be an
associative algebra; then we have a map
$$
q\colon H_*(\gl_\infty^{\Jac}(A);\C)\to H^{*+1}(\gl_\infty^{\fin}(A);\C).
$$
\end{example}

\subsection{}
The case we have discussed in Sect.~2.2 is the simplest case of the lifting
problem.
The next step is to extend our construction to cycles in
$H^1_{\Tr}(\A;\A)$. But we can do this only in the case $i=2$
(see~\S1). In the cases $i=4,6,8,\dots$ the most we can hope
is to solve the problem with some additional conditions. We thus assume
that the following conditions (i)--(iii) hold:

(i) $D_j\in\Der^{\Tr}_{\Ass}(\A)$ for all~$j$;

(ii) we solve problem \emph{only} for one monomial
$D_1\wedge\dots\wedge D_i$\ \ ($k=1$), and
$$
\text{$[D_a,D_b]=\ad Q_{ab}$ --- inner derivations ($Q_{ab}\in\A$);}
$$

(iii) $\Altl_{a,b,c}D_c(Q_{ab})=0$ for all $a$, $b$, $c$.

\begin{mainconjecture}[First Version]
Under the conditions \emph{(i)--(iii)} the lifting problem can be
solved for $\alpha^*_{D_1,\dots, D_i}$.
\end{mainconjecture}

We expect that
\begin{multline}
\Psi_{i+1}(A_1,\dots,A_{i+1})=\Altl_{A,D}\Tr(D_1A_1\cdot\ldots\cdot
D_iA_i\cdot A_{i+1})+\\
+(\text{terms, linear in $Q_{ab}$})
+(\text{terms, quadratic in $Q_{ab}$})
+\dots
\end{multline}
(by analogy with Theorem 1.1).

\begin{remark}
It follows easily from the Jacobi identity that $\Altl_{a,b,c}D_c(Q_{ab})$
lies in the center~$Z$ of Lie algebra~$\A$. Hence, codition (iii) is not
very strong; in particular, it holds for the Lie factoralgebra~$\A/Z$.
\end{remark}

Here the Conjecture will be proved only for $n=1,2$;
but we will conject in Sect.~4.6 an explicit formula for all~$n$.

The lifting formula for $n=2$ contains terms quadratic in~$Q_{ab}$.

\subsection{}
The main example when conditions (i)--(iii) hold is an
associative algebra $\Psi\Dif_n(S^1)$ of pseudo-differential
operators on~$(S^1)^n$ and $2\cdot n$ exterior derivations of
$
\Psi\Dif_n(S^1):\ln x_1,\dots,\ln x_n,\ln\partial_1,\dots,\ln\partial_n.
$
In fact, we can define derivation~$\ln\D$ (for
$\D\in\Psi\Dif_n(S^1)$) in much more generality. It seems to be
true that to do this we need only this condition:

\smallskip
\hangindent=2cm\hangafter=0\noindent
There exists  $\D^*\in\Psi\Dif_n(S^1)$ such that $[\D,\D^*]=1$.
\smallskip

But, first, we have to remember about conditions (i)--(iii),
and, second, we will see in Sect.~2.8 that in fact the choice
$\{\ln x_1,\dots,\ln x_n;\ln\partial_1,\dots,\ln\partial_n\}$ is, in
the some sense, most general.

\subsection{}
There is a deformation of the Lie algebra
$\Poiss_{2n}(S^1)$ to the Lie algebra
$\Psi\Dif_n(S^1)$. If we suppose that that the lifting problem
can be solved, it is interesting to find the Hamiltonian limit of
the lifting formulas. It is easy to see that in (10), this
limit depends only on the leading term, but not on the terms
containing~$Q_{ab}$.
Let us describe this cocycles.

Let $p_1$, \dots, $p_n$; $q_1$, \dots, $q_n$ be standard Poisson
coordinates, $\{p_i,q_j\}=\delta_{ij}$.

Let $\Tr\colon\Poiss_{2n}(S^1)\to\C$ be the linear functional
$\Tr=\Res_{p_1}\circ\ldots\circ\Res_{p_n}\circ\Res_{q_1}\circ\ldots\circ\Res_{
q_n}$.

\begin{lemma}
$\Tr\{f,g\}=0$ for  any $f,g\in\Poiss_{2n}(S^1)$.
\end{lemma}

\begin{proof}
$$
\{f,g\}=\sum_k\left(\nfracp f{p_k}\nfracp g{q_k}-\nfracp f{q_k}\nfracp
g{p_k}\right)=
\sum_k\left(\nfracp{}{p_k}\left(f\cdot\nfracp g{q_k}\right)-
\nfracp{}{q_k}\left(f\cdot\nfracp g{p_k}\right)\right).
$$
\end{proof}

\subsection{}
Any $\D\in\Poiss_{2n}(S^1)$, $\D\ne0$, determines an exterior
derivation $\ln\D\colon\Poiss_{2n}(S^1)\to\Poiss_{2n}(S^1)$ by
the formula
\begin{equation}
\{\ln\D,f\}=\nfrac{\{\D,f\}}{\D}.
\end{equation}

\begin{lemma}
For any $\D_1,\dots,\D_{2n}\in\Poiss_{2n}(S^1)$ the formula
$$
\Psi_{2n+1}^{\D_1,\dots,\D_n}(f_1,\dots,f_{2n+1})=
\Altl_f\Tr(\{\ln\D_1,f_1\}\cdot\ldots\cdot\{\ln\D_{2n},f_{2n}\}\cdot f_{2n+1})
$$
defines a $(2n+1)$-cocycle on $\Poiss_{2n}(S^1)$.
\end{lemma}

\begin{proof}
We will proof, that
$$
\Psi_{2n+1}^F(f_1,\dots,f_{2n+1})=\Altl_f\Tr(F\cdot\nfracp{f_1}{p_1}\cdot
\ldots\cdot\nfracp{f_n}{p_n}\cdot\nfracp{f_{n+1}}{q_1}\cdot\ldots\cdot
\nfracp{f_{2n}}{q_n}\cdot f_{2n+1})
$$
defines a $(2n+1)$-cocycle on $\Poiss_{2n}(S^1)$.
It is easy to see, that
\begin{gather*}
\Psi_{2n+1}^{\D_1,\dots,\D_{2n}}=\Psi^F_{2n+1}\\
\text{for}\quad F=\nfrac{\det\left(\nfracp{\D_i}{\xi_j}\right)}{\D_1\cdot\ldots\cdot\D_{2n}}\quad
\text{(here $\xi_1=p_n$, \dots, $\xi_n=p_n$, $\xi_{n+1}=q_1$, \dots,
$\xi_{2n}=q_n$)}.
\end{gather*}
We have:
\begin{multline}
(d\Psi^F_{2n+1})(f_{2n+2},f_1,\dots,f_{2n+1})=\\
=\nfrac12\Altl_f\Tr\left[F\cdot\nfracp{}{p_1}\{f_{2n+2}f_1\}\cdot
\nfracp{}{p_2}f_2\cdot\nfracp{}{p_3}f_3\cdot\ldots
+F\cdot\nfracp{}{p_1}f_1\cdot\nfracp{}{p_2}\{f_{2n+2},f_2\}\cdot\nfracp{}{p_3}
f_3\cdot\right.\\
\cdot\ldots
+\ldots+
\left.+F\cdot\nfracp{}{p_1}f_1\cdot\ldots\cdot\nfracp{}{q_n}f_{2n}\cdot\{f_{2n
+2},f_{2n+1}\}\right].
\end{multline}
We subtract from the r.h.s. of (12) the expression
$$
\nfrac12\Altl_f\Tr\left(\left\{f_{2n+2},
F\cdot\nfracp{}{p_1}f_1\cdot\ldots\cdot\nfracp{}{q_n}f_{2n}\cdot
f_{2n+1}\right\}\right),
$$
which is equal to zero by Lemma 2.6.

We have:
\begin{multline}
\text{r.h.s. (12)}=-\nfrac12\Altl_f
\Tr\left(\{f_{2n+2},
F\}\cdot\nfracp{}{p_1}f_1\cdot\ldots\cdot\nfracp{}{q_n}f_{2n}\cdot
f_{2n+1}\right)\\
+\nfrac12\Altl_f\Tr\left[F\cdot\{\nfracp{}{p_1}f_{2n+2},f_1\}\cdot
\nfracp{}{p_2}f_2\cdot\ldots
+F\cdot\nfracp{}{p_1}f_1\cdot\nfracp{}{p_2}\{f_{2n+2},f_2\}\cdot\ldots+\ldots
\right].
\end{multline}
The second summand is equal to the first one by Lemma 2.5 and the
symmetry, and we have:
\begin{equation}
(d\Psi^F_{2n+1})(f_{2n+2},f_1,\dots,f_{2n+1})=
\Altl_f\Tr\left(\{F,f_{2n+2}\}\cdot\nfracp{}{p_1}f_1\cdot\ldots
\nfracp{}{q_n}f_{2n}\cdot f_{2n+1}\right).
\end{equation}
Furthermore $\{F,f_{2n+2}\}=\suml_k\left(\nfracp F{p_k}\nfracp{f_{2n+2}}{q_k}-
\nfracp F{q_k}\cdot\nfracp{f_{2n+2}}{p_k}\right)$ and every
summand in (14) is equal to zero via the alternation.

\subsection{}
Although Lemma 2.7 gives us a lot of cocycles, it is easy to
see that $\Psi_{2n+1}^{F_1}\sim\Psi_{2n+1}^{F_2}$ when $\Res
F_1=\Res F_2$; therefore, we have only one interesting cocycle,
which is $\Psi_{2n+1}^{\ln p_1,\dots,\ln p_n,\ln q_1,\dots, ln
q_n}$. We will denote it by $\Psi_{2n+1}^0$.

\begin{lemma}
\begin{enumerate}
\item If the lifting problem for $\Psi\Dif_n(S^1)$ can be
solved, the Hamiltonian limit of any lifting is~$\Psi_{2n+1}^0$;
\item $\Psi_{2n+1}^0(1\wedge p_1\wedge\dots\wedge p_n\wedge
q_1\wedge\dots\wedge q_n)\ne0$.
\end{enumerate}
\end{lemma}

\end{proof}

This Lemma implies that \emph{the lifting problem \emph= the deformation problem.}

On the other hand, we see that we don't need try to generalize
the construction for $\Psi\Dif_n(S^1)$ to other derivations
of the form~$\ln\D$; the case $D_1=\ln x_1$, \dots,
$D_{2n}=\ln\partial_n$ is most general.

\begin{corollary}
Any $(2n+1)$-cocycle on the Lie algebra $\Psi\Dif_n(S^1)$ of the
form~\emph{(10) ($D_1=\ln x_1$, \dots, $D_{2n}=\ln\partial_n$)} is
not cohomologous to zero.
\end{corollary}

\epr

Another proof follows from the fact that every cocycle of the
form~(10) is a lifting of $\alpha^*_{\ln x_1,\dots,\ln\partial_n}$,
which is a nonzero element in $H^{2n}(\Psi\Dif_n(S^1);\Psi\Dif_n(S^1)^*)$.

\section{Computation for $n=2$}

\subsection{}
The main result of this Section is the following

\begin{theorem}
Let $\A$ be an associative algebra, $\Tr\colon\A\to\C$ be a
trace on~$\A$ and $D_1$, $D_2$, $D_3$, $D_4$ --- four
\emph(exterior\emph) derivations on~$\A$, which satisfy
conditions \emph{(i)--(iii)} from Sect.~\emph{2.4}. Then
\begin{multline}
\Psi_5(A_1,A_2,A_3,A_4,A_5)=
\Altl_{A,D}\{
D_1A_1\cdot D_2A_2\cdot D_3A_3\cdot D_4A_4\cdot A_5 \\
{\left.{\begin{aligned}
&+A_1Q_{12}A_2\cdot D_3A_3\cdot D_4A_4\cdot A_5\\
&+D_1A_1\cdot A_2\cdot Q_{23}\cdot A_3\cdot D_4A_4\cdot A_5\\
&+D_1A_1\cdot D_2A_2\cdot A_3\cdot Q_{34}\cdot A_4\cdot A_5
\end{aligned}}\right]}\
\text{terms, linear in $Q_{ij}$}\\
+A_1Q_{12}A_2A_3Q_{34}A_4A_5\}]\ \text{term, quadric in $Q_{ij}$}
\end{multline}
is a $5$-cocycle on the Lie algebra~$\A$, which is the lifting
of $\alpha^*_{D_1,D_2,D_3,D_4}$ \emph(see Sect.~\emph{2.1)}. In the
case $\A=\Psi\Dif_2(S^1)$,\ \ $D_1=\ln\partial_1$, $D_2=\ln
x_1$, $D_3=\ln\partial_2$, $D_4=\ln x_2$,\ \ $\Psi_5$~is not
cohomologous to zero\emph: $\Psi_5(1\wedge x_1\wedge
x_2\wedge\partial_1\wedge\partial_2)\ne0$.
\end{theorem}

\begin{remark} $\ad Q_{ij}=[D_i,D_j]$; however, in formula (15),
we don't alternate symbols~$i$ and~$j$ in~$Q_{ij}$.
\end{remark}

We will give a sketch of the proof of this Theorem in Sect.\ 3.2--3.5.

\subsection{}

\begin{lemma}
\begin{multline*}
(\mathrm i)\ d(\Altl_{A,D}\Tr(D_1A_1\cdot A_2\cdot A_3\cdot D_2A_4\cdot Q\cdot
A_5)(A_6,A_1,\dots,A_5)=\\
=\Altl_{A,D}\Tr(D_1A_1\cdot A_2\cdot A_3\cdot D_2A_6\cdot A_4\cdot Q\cdot
A_5)\\
\end{multline*}
\begin{multline*}
(\mathrm{ii})\ d(\Altl_{A,D}\Tr(D_1A_1\cdot D_2A_2\cdot A_3\cdot A_4\cdot Q\cdot
A_5)(A_6,A_1,\dots,A_5)=\\
=\Altl_{A,D}\Tr(D_1A_1\cdot D_2A_6\cdot A_2\cdot A_3\cdot A_4\cdot
Q\cdot A_5)\\
\end{multline*}
\begin{multline*}
(\mathrm{iii})\ d(\Altl_{A,D}\Tr(D_1A_1\cdot D_2A_2\cdot A_3\cdot A_4\cdot Q\cdot
A_5)(A_6,A_1,\dots,A_5)=\\
=\Altl_{A,D}\Tr(D_1A_1\cdot D_2A_2\cdot A_3\cdot Q\cdot A_6 \cdot A_4\cdot
A_5)
\end{multline*}
\end{lemma}

\begin{proof}
Straightforward (see Sect.~1.3).
\end{proof}

\subsection{}

\begin{lemma}
\begin{multline*}
(\mathrm i)\ d(\Altl_{A,D}\Tr(D_1A_1\cdot A_2\cdot A_3\cdot A_4\cdot D_2A_5\cdot
Q))(A_6,A_1,\dots,A_5)=\\
=\Altl_{A,D}\Tr(D_1A_6\cdot A_1\cdot A_2\cdot A_3\cdot A_4\cdot
D_2A_5\cdot Q-D_1A_1\cdot A_2\cdot A_3\cdot A_4\cdot
D_2A_5\cdot[A_6,Q])\\
\end{multline*}
\begin{multline*}
(\mathrm{ii})\ d(\Altl_{A,D}\Tr(D_1A_1\cdot A_2\cdot D_2A_3\cdot A_4\cdot A_5\cdot
Q))(A_6,A_1,\dots,A_5)=\\
=\Altl_{A,D}\Tr(D_1A_6\cdot A_1\cdot A_2\cdot D_2A_3\cdot A_4\cdot
A_5\cdot Q-D_1A_1\cdot A_2\cdot D_2A_3\cdot A_4\cdot
A_5\cdot[A_6,Q])\\
\end{multline*}
\begin{multline*}
(\mathrm{iii})\ d(\Altl_{A,D}\Tr(D_1A_1\cdot A_2\cdot D_2A_3\cdot Q\cdot A_4\cdot
A_5))(A_6,A_1,\dots,A_5)=\\
=\Altl_{A,D}\Tr(D_1A_6\cdot A_1\cdot A_2\cdot D_2A_3\cdot Q\cdot A_4\cdot
A_5-D_1A_1\cdot A_2\cdot D_2A_3\cdot[A_6,Q]\cdot A_4\cdot
A_5)
\end{multline*}
\end{lemma}

\epr

\subsection{}

\begin{lemma}
\begin{multline*}
(\mathrm i)\ d(\Altl_{A,D}\Tr(D_1A_1\cdot D_2A_2\cdot D_3A_3\cdot D_4A_4\cdot
A_5))(A_6,A_1,\dots,A_5)=\\
=\Altl_{A,D}\Tr(D_1A_1\cdot D_2A_2\cdot D_3A_3\cdot D_4A_4\cdot
A_5\cdot A_6)\\
\end{multline*}
\begin{multline*}
(\mathrm{ii})\ d(\Altl_{A,D}\Tr(Q_1\cdot Q_2\cdot A_1\cdot A_2\cdot A_3\cdot A_4\cdot
A_5))(A_6,A_1,\dots,A_5)=\\
=-\Altl_A\Tr(Q_1\cdot Q_2\cdot A_1\cdot  A_2\cdot A_3\cdot A_4\cdot
A_5\cdot A_6)\\
\end{multline*}
\begin{multline*}
(\mathrm{iii})\ d(\Altl_{A,D}\Tr(Q_1\cdot A_1\cdot Q_2\cdot A_2\cdot A_3\cdot A_4\cdot
A_5))(A_6,A_1,\dots,A_5)=\\
=\Altl_A\Tr(Q_1\cdot A_6\cdot A_1\cdot Q_2\cdot A_2\cdot A_3\cdot
A_4\cdot A_5)
\end{multline*}
\end{lemma}

\epr

\subsection{}

To prove Theorem 3.1, we consider an expression in $A_1$,
\dots,~$A_6$ which a priori is equal to zero (we use
$\Tr(DA)=0$ --- condition~(i) from Sect.~2.4), and then prove
that this expression is a coboundary of $\Psi_5(A_1,\dots,A_5)$.
Then $\Psi_5$ is a $5$-cocycle. Our main tools are Lemmas 3.2--3.4.

Suppose that
\begin{multline*}
E_1=\Altl_{A,D}\Tr(2D_2(Q_{43}\cdot D_1A_1\cdot A_2\cdot A_3\cdot A_4\cdot
A_5\cdot A_6)+\\
+2D_2(D_1A_1\cdot Q_{43}\cdot A_2\cdot A_3\cdot A_4\cdot A_5\cdot A_6))\\
\end{multline*}
\begin{multline*}
E_2=\Altl_{A,D}\Tr(D_2(D_1A_1\cdot  A_2\cdot A_3\cdot Q_{43}\cdot A_4\cdot
A_5\cdot A_6)+\\
+D_2(D_1A_1\cdot  A_2\cdot A_3\cdot A_4\cdot Q_{43}\cdot A_5\cdot A_6))\\
\end{multline*}
$$
E_3=\Altl_{A,D}\Tr(D_4(D_1A_1\cdot A_2\cdot A_3\cdot A_4\cdot Q_{43}\cdot
A_5\cdot A_6))
$$

Note that $E_1=E_2=E_3=0$.

Furthermore suppose that
\begin{gather*}
I=\Altl_{A,D}\Tr(D_2((D_1A_1\cdot A_2\cdot A_3\cdot A_4)\cdot D_4(A_6\cdot
D_3A_5)))\\
II=\Altl_{A,D}\Tr(D_2((A_3\cdot A_4\cdot D_1A_1\cdot A_2)\cdot D_4(A_6\cdot
D_3A_5)))\\
III=\Altl_{A,D}\Tr(D_2((A_2\cdot A_3\cdot A_4\cdot D_1A_1)\cdot
D_4(D_3A_5\cdot A_6)))\\
IV=\Altl_{A,D}\Tr(D_2((A_4\cdot D_1A_1\cdot A_2\cdot A_3)\cdot D_4(D_3A_5\cdot
A_6)))\\
\end{gather*}

Note also that $I=II=III=IV=0$.

We state that
\begin{equation}
E_1-E_2+\nfrac43E_3-I+II-III+IV=d\wt\Psi_5
\end{equation}
where
\begin{multline*}
\wt\Psi_5
=\Altl_{A,D}\Tr(
-2D_1A_1\cdot D_2A\cdot D_3A_3\cdot D_4A_4\cdot A_5]\ \text{term without
$Q_{ij}$}\\
{\left.{\begin{aligned}
&+D_1A_1\cdot D_2A_2\cdot Q_{43}\cdot A_3\cdot A_4\cdot A_5\\
&+D_1A_1\cdot D_2A_2\cdot A_3\cdot A_4\cdot A_5\cdot Q_{43}\\
&+2\cdot D_1A_1\cdot Q_{43}\cdot D_2A_2\cdot A_3\cdot A_4\cdot A_5\\
&-D_1A_1\cdot A_2\cdot D_2A_3\cdot Q_{43}\cdot A_4\cdot A_5\\
&-D_1A_1\cdot A_2\cdot D_2A_3\cdot A_4\cdot A_5\cdot Q_{43}
\end{aligned}}\right]}\ \text{terms, linear in $Q_{ij}$}\\
{\left.{\begin{aligned}
&-4Q_{21}\cdot Q_{43}\cdot A_1\cdot A_2\cdot A_3\cdot A_4\cdot A_5\\
&+2Q_{21}\cdot A_1\cdot A_2\cdot Q_{43}\cdot A_3\cdot A_4\cdot A_5\\
&-2[Q_{21},A_1]\cdot[Q_{43},A_2]\cdot A_3\cdot A_4\cdot A_5
\end{aligned}}\right]}\ \text{terms, quadratic in $Q_{ij}$}
\end{multline*}
To prove this, we need condition (iii) from Sect.~2.4:
$$
\Altl_{i,j,k}D_k(Q_{ij})=0\quad\text{for all $i$, $j$, $k$}.
$$
The proof of (16) is a very long straightforward calculation,
using Lemmas 3.2--3.4 and other observations.

Furthermore, after simple manipulations
$$
\wt\Psi_5-\Altl_{A,D}\Tr(D_2(D_1A_1\cdot Q_{43}\cdot A_2\cdot  A_3\cdot
A_4\cdot A_5))-
\Altl_{A,D}\Tr(D_1(A_1\cdot Q_{43}\cdot D_2A_2\cdot A_3\cdot A_4\cdot A_5))
$$
will have a form as in Th. 3.1.
\epr

\section{Lifting formulas: the general case}

\subsection{}
In this Section the lifting formulas for any number of derivations
$D_1$, \dots,~$D_l$ (with conditions (i)--(iii) from Sect.~2.4)
appear. The fact that polylinear skew-symmetric functions
on~$\A$, defined by these formulas, are \emph{cocycles}, is the
Main Conjecture of this Section. As we have seen in~\S3, in our
situation it is highly nontrivial to check that a given formula
in fact defines a cocycle.

The main idea is the following. We consider the expression
\begin{equation}
\Altl_{A,D}\Tr(D_1A_1\cdot D_2A_2\cdot\dots\cdot D_{2n}A_{2n}\cdot A_{2n+1})
\end{equation}
where $D_1$, \dots,$D_{2n}$ are \emph{inner} derivations,
$D_iA=D_i\cdot A-A\cdot D_i$. We want to add to (17) some terms
containing $Q_{ij}=[D_i,D_j]$ and in this way obtain a cocycle. It
is meaningless from the cohomological viewpoint, because this
cocycle will be cohomologous to~$0$ (see Remark~4.2 and
Remark~4.5), but using this trick we obtain explicit formulas.

We have:
$$
(17)=\Altl_{A,D}\Tr((D_1\cdot A_1-A_1\cdot D_1)\cdot(A_2\cdot D_2-D_2\cdot
A_2)\cdot\ldots\cdot
(D_{2n}\cdot A_{2n}-A_{2n}\cdot D_{2n})\cdot A_{2n+1})
$$

Note also that our formulas are true for any number of
derivations, not only for an even number.

\subsection{}
\begin{lemma}
Let $\A$ be an associative algebra; $D_1, \dots, D_{2n}\in\A$. Then
$$
\alpha(A_1,\dots,A_{2n+1})=\Altl_{A,D}\Tr(D_1\cdot A_1\cdot
D_2\cdot A_2\cdot\ldots\cdot D_{2n}\cdot A_{2n}\cdot A_{2n+1})
$$
is a $(2n+1)$-cocycle on the Lie algebra~$\A$.
\end{lemma}

\begin{proof}
\begin{multline}
d\alpha(A_{2n+2},A_1,\dots,A_{2n+1})=\\
=\nfrac14\Altl_{A,D}\Tr\{
D_1\cdot[A_{2n+2},A_1]\cdot D_2\cdot A_2\cdot{\ldots}\cdot D_{2n}
\cdot[A_{2n},A_{2n+1}]\\
{\begin{aligned}
&+D_1\cdot A_1\cdot D_2 \cdot[A_{2n+2},A_2]\cdot D_3\cdot
A_3\cdot{\ldots}\cdot D_{2n}\cdot[A_{2n},A_{2n+1}]\\
&+{\ldots}\\
&+D_1\cdot A_1\cdot D_2\cdot A_2\cdot{\ldots}\cdot
D_{2n-1}\cdot[A_{2n+2},A_{2n-1}]\cdot D_{2n}\cdot[A_{2n},A_{2n+1}]\\
&+D_1\cdot A_1\cdot{\ldots}\cdot D_{2n}\cdot[[A_{2n+2},[A_{2n},A_{2n+1}]]\}.
\end{aligned}}
\end{multline}
The last summand in (18) is equal to~$0$ by the Jacobi identity.
The remaining summands cancel by the trace property and
symmetry conditions. For example,
\begin{multline*}
\Altl_{A,D}\Tr\{D_1\cdot[A_{2n+2},A_1]\cdot D_2\cdot A_2\cdot\ldots\cdot
D_{2n}\cdot[A_{2n},A_{2n+1}]+\\
+D_1\cdot A_1\cdot\ldots\cdot D_{2n-1}\cdot[A_{2n+2},A_{2n-1}]\cdot
D_{2n}\cdot[A_{2n},A_{2n+1}]\}=0.
\end{multline*}
Note also that the proof remains the same in the case when there are an odd number
of~$D_i$.
\end{proof}

\begin{remark}
Indeed, $\alpha$ is a coboundary, $\alpha=a\cdot d(D_1\cdot
A_1\cdot\ldots\cdot D_{2n}\cdot A_{2n})$\ \ ($a\in\Z$).
\end{remark}

\subsection{}
Let $D_i\in\A$ for all~$i$. Then
\begin{equation}
\Altl_{A,D}\Tr([D_1,A_1]\cdot\ldots\cdot[D_{2n},A_{2n}]\cdot A_{2n+1})=
k\alpha+S\quad(k\in\Z)
\end{equation}
where
\begin{equation}
S=(\text{sum of the terms which contains $D_i\cdot D_{i+1}$ for some~$i$})
\end{equation}
We replace in the all summands of~$S$\ \ $D_iD_{i+1}$ by $[D_i,D_{i+1}]$
because of the symmetry condition. Hence, these terms have the same form
as the terms in the lifting formulas, and our problem is to
represent $S$ as a sum of the terms of the form
\begin{multline*}
\Altl_{A,D}\Tr([D_1,A_1]\cdot\ldots\cdot A_{i_1}\cdot
Q_{i_1,i_1+1}A_{i_1+1}\cdot[D_{i_1+2},A_{i_1+2}]\\
\cdot\ldots\cdot
A_{i_2}\cdot Q_{i_2,i_2+1}\cdot 
A_{i_2+1}\cdot[D_{i_2+2},A_{i_2+2}]\cdot\ldots).
\end{multline*}
Then, when we subtract $S$ from (19) we obtain~$k\alpha$, which
is a cocycle by Lemma~2.1. Indeed, in this case ($D_i\in A$),\ \
$k\alpha$~is a coboundary. But it turns out that when the $D_i$ are
\emph{exterior} derivations \emph{satisfying conditions}~(i)--(iii)
from Sect.~2.4 and $n=1,2$, the lifting formulas from \S1,3 have the
same form. This fact allows us to formulate the Main Conjecture
in Sect.~4.6.

\subsection{}
We will denote summands from~$S$ by closed intervals with
marked points. If $2n$ is the total number of~$D_i$, the length of
the interval is $2n-2$, and some integral points on it are
marked. Let us denote by $1$, \dots, $2n-1$ the integral points of the
interval; the point~$i$ is marked, if in the corresponding
summand in~$S$,\ \ $D_i$ and $D_{i+1}$ are neighbors (i.e.
not separated by $A_i$ or $A_{i+1}$, or $A_i\cdot A_{i+1}$).
For example, the interval\newline
\smallskip
\centerline{$\vcenter{\epsfbox{pict.1}}$}
\smallskip

\noindent
corresponds to the expression
$$
\Altl_{A,D}\Tr(A_1\cdot D_1\cdot D_2\cdot A_2\cdot[D_3,A_3]\cdot[D_4,A_4]\cdot A_5).
$$
It is clear that the distance between the marked points
is~$\ge2$. In the general case, the interval\newline
\smallskip
\centerline{$\vcenter{\epsfbox{pict.2}}$}
\smallskip

\noindent
corresponds to the expression
\begin{multline}
\Altl_{A,D}\Tr([D_1,A_1]\cdot[D_2,A_2]\cdot\ldots\cdot A_{i_1}\cdot
D_{i_1}\cdot D_{i_1+1}\cdot A_{i_1+1}\cdot
[D_{i_1+2},A_{i_1+2}]\\
\cdot\ldots\cdot A_{i_2}\cdot D_{i_2}\cdot
D_{i_2+1}A_{i_2+1}\cdot\ldots)
\end{multline}
Let $N$ be the total number of the marked points, $1\le N\le n$.

The first summand in $S$ is a sum of all intervals with $N=1$, with sign 
``$-$'' ($2n-1$ intervals):
\begin{figure}[h]
$$
S_1=-\left(\vcenter{\epsfbox{pict.3}}\right)
$$
\caption{$n=3$}
\end{figure}

\noindent
Some summands of $S$ are counted in~$S_1$ with multiplicities~$>1$. We
subtract from~$S_1$ all summands with $N=2$ pairs of~$D_i$, and
then we add every interval with $N=2$ marked points with its
sign in~$S$. We have:
\begin{figure}[h]
$$
S_2=-2\cdot\left(\vcenter{\epsfbox{pict.4}}\right)+1\cdot\left(\vcenter{\epsfbox{pict.4}}\right)
$$
\caption{$S_2$ for $n=3$}
\end{figure}

\noindent
We wrote the coefficient $-2$ because every interval with $N=2$
marked points is ``contained'' in two intervals with $N=1$ marked
points with the sign ``$-$''; for example, the interval $\vcenter{\epsfbox{pict.5}}$
``contains'' in $\vcenter{\epsfbox{pict.6}}$ and in $\vcenter{\epsfbox{pict.7}}$. In the sum $S_1+S_2$
only the terms with $N\ge3$ pairs of
$D_i$
are
counted incorrectly.
In our example ($n=3$) the maximal value for $N$ is
$N=3$, and there exists just one unique interval with $N=3$:
$\vcenter{\epsfbox{pict.8}}$.

\pagebreak

It contains with the sign ``$+$'' in $3$ intervals with $N=1$ marked points, 
namely, in
\begin{figure}[h]
\centerline{$\vcenter{\epsfbox{pict.9}}$}
\caption{}
\end{figure}

\noindent
and with the sign ``$-$'' in $3$ intervals with $N=2$ marked points:
\begin{figure}[h]
\centerline{$\vcenter{\epsfbox{pict.10}}$}
\caption{}
\end{figure}

\noindent
Hence,
\def\fignja{\vcenter{\epsfbox{pict.8}}}
$$
S_3=3\cdot\fignja-3\cdot\fignja-1\cdot\fignja=-1\cdot\fignja.
$$

Note that all coefficients in $S_1$, $S_2$, $S_3$ are equal to~$-1$.

We have: $S=S_1+S_2+S_3$.

\subsection{}
Let us denote by $\Sigma_{N,2n-1}$ the sum of \emph{all} intervals of
length $2n-2$ with $N$ marked points ($1\le N\le n$, the
distance between marked points is $\ge2$).

\begin{lemma}
$$
S=-\Sigma_{1,2n-1}-\Sigma_{2,2n-1}-\ldots-\Sigma_{n,2n-1}.
$$
\end{lemma}

\begin{proof}
Suppose that the summands in $S$ with $N\le k$ pairs of~$D_i$
are counted correctly in $S^{(k)}=-\Sigma_{1,2n-1}-\ldots-\Sigma_{k,2n-1}$.
Then every summand with $N=k+1$ pairs of~$D_i$ is ``contained''
in~$\Sigma_i$ with the multiplicity~$C_{k+1}^i$ and with the
sign~$(-1)^{i-1}$\ \ $i\le k$.

Hence we add to~$S^{(k)}$
\begin{equation}
(-1)\cdot\left(\sum_{i=1}^k(-1)^iC_{k+1}^i\right)\Sigma_{k+1,2n-1}
+(-1)^{k+1}\Sigma_{k+1,2n-1}
\end{equation}
It is clear that in the sum $S^{(k)}+\text{(22)}$ all terms with
$N\le k+1$ pairs of~$D_i$ are correctly calculated; on the other
hand, $\text{(22)}=-1\Sigma_{k+1,2n-1}$ by the binomial formula
for \hbox{$(1-1)^{k+1}=0$}.
\end{proof}

\begin{corollary}
For $D_1,\dots,D_{2n}\in\A$,
$$
\Altl_{A,D}\Tr([D_1,A_1]\cdot\ldots\cdot[D_{2n},A_{2n}]\cdot A_{2n+1})
+\Sigma_{1,2n-1}+\ldots+\Sigma_{n,2n-1}
$$
is a cocycle on the Lie algebra~$\A$.
\end{corollary}

\begin{proof}
Follows from Lemma~4.2.
\end{proof}

\begin{remark}
It follows from Remark 4.2 that this cocycle is cohomologous to~$0$.
\end{remark}

\subsection{}
\begin{mainconjecture}[Second Version]
Let $\A$ be an associative algebra with~$\Tr$,\ \ $D_1$, \dots,
$D_{2n}$ be its derivations satisfying conditions \emph{(i)--(iii)}
from Sect.~\emph{2.4}. Then \emph(in the notations of Sect. \emph{4.4,~4.5)}
\begin{multline*}
\Psi_{2n+1}(A_1,\dots,A_{2n+1})=\\
=\Altl_{A,D}\Tr(D_1A_1\cdot\ldots\cdot D_{2n}A_{2n}\cdot A_{2n+1})
+\Sigma_{1,2n-1}+\Sigma_{2,2n-1}+\ldots+\Sigma_{n,2n-1}
\end{multline*}
is a cocycle on the Lie algebra~$\A$ \emph{(and, hence, it solves
the lifting problem)}.
\end{mainconjecture}

\begin{remark}
Let $\A$ be an associative algebra, $\h$~be a two-sided ideal
in~$\A$,\ \ $\Tr$ --- a trace on~$\h$,\ \
 $\A$~acts on~$\h$ by the adjoint action (see Ex.
2.3.2). Then the above formula is true without conditions
(ii)--(iii) from Sect.~2.4, but any  cocycle obtained in
this way is cohomologous to~$0$.
\end{remark}

\subsection{}
\begin{example}[$n=3$]
In \S1 we obtained the lifting formula for $\Dif_1$, in \S3 we found it
for~$\Dif_2$. Main Conjecture~4.6 states that the following
formula is a lifting formula for $\Dif_3$:
\begin{multline*}
\Psi_3(A_1,\dots,A_7)=\Altl_{A,D}\Tr\{D_1A_1\cdot\ldots\cdot D_6A_6\cdot A_7\\
\left.
{\begin{aligned}
&+A_1\cdot Q_{12}\cdot A_2\cdot D_3A_3\cdot D_4A_4\cdot D_5A_5\cdot 
D_6A_6\cdot A_7\\
&+D_1A_1\cdot A_2\cdot Q_{23}\cdot A_3\cdot D_4A_4\cdot D_5A_5\cdot 
D_6A_6\cdot A_7\\
&+D_1A_1\cdot D_2A_2\cdot A_3\cdot Q_{34}\cdot A_4\cdot D_5A_5\cdot 
D_6A_6\cdot A_7\\
&+D_1A_1\cdot D_2A_2\cdot D_3A_3\cdot A_4\cdot Q_{45}\cdot A_5\cdot 
D_6A_6\cdot A_7\\
&+D_1A_1\cdot D_2A_2\cdot D_3A_3\cdot D_3A_4\cdot A_5\cdot Q_{56}\cdot 
A_6\cdot A_7\\
\end{aligned}}
\right]\ \text{terms, linear in $Q_{ij}$}
\\
\left.
{\begin{aligned}
&+A_1\cdot Q_{12}\cdot A_2\cdot A_3\cdot Q_{34}\cdot A_4\cdot D_5A_5\cdot 
D_6A_6\cdot A_7\\
&+A_1\cdot Q_{12}\cdot A_2\cdot D_3A_3\cdot A_4\cdot Q_{45}\cdot A_5\cdot 
D_6A_6\cdot A_7\\
&+A_1\cdot Q_{12}\cdot A_2\cdot D_3A_3\cdot D_4A_4\cdot A_5\cdot Q_{56}\cdot 
A_6\cdot A_7\\
&+D_1A_1\cdot A_2\cdot Q_{23}\cdot A_3\cdot A_4\cdot Q_{45}\cdot A_5\cdot
D_6A_6\cdot A_7\\
&+D_1A_1\cdot A_2\cdot Q_{23}\cdot A_3\cdot D_4A_4\cdot A_5\cdot Q_{56}\cdot
A_6\cdot A_7\\
&+D_1A_1\cdot D_2A_2\cdot A_3\cdot Q_{34}\cdot A_4\cdot A_5\cdot Q_{56}\cdot 
A_6\cdot A_7\\
\end{aligned}}
\right]\ \text{terms, quadratic in $Q_{ij}$}
\\
+A_1\cdot Q_{12}\cdot A_2\cdot A_3\cdot Q_{34}\cdot A_4\cdot A_5\cdot 
Q_{56}\cdot A_6\cdot A_7\}
]\ \text{term, cubic in $Q_{ij}$}
\end{multline*}
Recall that we don't alternate symbols $i$, $j$ in $Q_{ij}$.
\end{example}

\end{document}